\documentclass[12pt,a4paper]{article}
\usepackage{lineno,hyperref}
\usepackage[latin5]{inputenc}
\usepackage{amsmath}
\usepackage{graphicx}
\usepackage{subfigure}
\usepackage{caption}
\usepackage{float}
\usepackage{mathtools}
\usepackage{color}
\usepackage{cite}
\usepackage{mathrsfs}
\usepackage[]{geometry}
\usepackage{subfigure}
\DeclareMathOperator{\sech}{sech}
\title{Extended B-Spline Collocation Method For KdV-Burgers Equation}
\author{Ozlem Ersoy$^{a}$, Alper Korkmaz$^{b,}$\thanks{akorkmaz@karatekin.edu.tr}, Idiris Dag$^{c}$  \\
$^{a}${\scriptsize Department of Mathematics \& Computer, Eskisehir Osmangazi University, 26480, Eskisehir, Turkey.} \\
$^{b}${\scriptsize Department of Mathematics, Çankırı Karatekin University, 18200, Çankırı, Turkey.}\\
$^{c}${\scriptsize Department of Computer Engineering, Eskisehir Osmangazi University, 26480, Eskisehir, Turkey.}}
\begin{document}
\maketitle
\begin{abstract}

\noindent
The extended form of the classical polynomial cubic B-spline function is used to set up a collocation method for some initial boundary value problems derived for the Korteweg-de Vries-Burgers equation. Having nonexistence of third order derivatives of the cubic B-splines forces us to reduce the order of the term $u_{xxx}$ to give a coupled system of equations. The space discretization of this system is accomplished by the collocation method following the time discretization with Crank-Nicolson method. Two initial boundary value problems, one having analytical solution and the other is set up with a non analytical initial condition, have been simulated by the proposed method.

\end{abstract}




\section{Introduction}
\noindent
Consider the Korteweg-de Vries - Burgers (KdVB) equation of the form
\begin{equation}
u_{t}+\varepsilon uu_{x}-\vartheta u_{xx}+\mu u_{xxx}=0 \ a\leq x\leq b
\label{kdvb}
\end{equation}
where $u=u(x,t)$, subscripts denote the derivatives and $\varepsilon$,$\nu$ and $\mu$ are real coefficients with the property $\varepsilon\nu\mu \neq 0$\cite{johnson1,karpman1,kalinowski1}.
\noindent
The KdVB equation for dispersive and dissociative media is derived in various cases such as existence and absence of a complete system of eigenvector or for waves in plasma with Hall dispersion and Joule dissipation by Ruderman\cite{ruderman1}. Including both viscosity and inertia terms simultaneously causes long gravity waves to be governed by the KdVB equation. This result was obtained by deriving the balance equations for an incompressible viscous fluid starting from the column model approximation in the study of Bampi and Morro\cite{bampi1}. 

\noindent
The Weierstrass $\mathscr{P}$-function solution, which can be expressed in terms of Jacobi cosine functions, for the KdVB \ref{kdvb} was found by Kalinowski and Grundland \cite{kalinowski1}. In the study, the Riemann invariants for non-homogenous systems of the first order partial differential equations method was implemented to obtain that solution. Brugarino and Pantano \cite{burigarino1} obtaned the Jacobi elliptic type solution of the nonhomogenous KdVB equation with variable coefficients. Some travelling solitary wave solutions to compound KdVB equation was developed by using automated method based on an ansatz of a polynomial in terms of a $\tanh$ function\cite{parkes1}. A particular analytical solution for the KdVB equation was also obtained by using variable transformations and proofs of theorems\cite{shu1}. The solutions in the stationary waves form were determined for the generalized KdVB equation containing nonconservative terms of linear pump, linear HF dissipation, and nonlinear dissipation\cite{gromov1}. Malkov \cite{malkov1} constructed the asymptotic travelling wave solution, representing a shock-train, of the KdVB equation driven by the long scale periodic driver. The stability of travelling wave solutions of the KdVB equation was studied for various values of the viscosity coefficient $\alpha$\cite{rego1}.

\noindent
The exact solutions consisting of powers of some hyperbolic functions of the KdVB equation were derived by homogenous balance technique in \cite{wang1}. These solutions are the combination of a belly-shape and kink-type solitary waves. Zhao \cite{zhao1} implemented the hyperbolic function method and the Wu elimination technique for the new type solitary wave and periodic solutions containing some blow-up solutions of the KdVB. After reducing the KdVB equation to an ordinary differential equation by the compatible wave transformation,  he determined the coefficients and parameters in the predicted solution by substituting it into this equation. Yuanxi and Jiashi\cite{yuanxi} developed many solitary wave solutions for the KdVB equation by the superposition method based on the analysis on the features of the Burgers, the KdV and the KdV-Burgers equations. A generalized $\tanh$ function method based on Riccati equation was derived to obtain multiple soliton solutions containing some trigonometric, hyperbolic and complex functions for the KdVB equation\cite{wazzan1}. Soliman\cite{soliman1} obtained the $\tanh$-type, $\coth$-type exact solutions of the equation by the modified extended $\tanh$ method.  $\tanh$-type solution for the KdVB equation was obtained by using $G'/G$ expansion method\cite{wang2}. Some compact and non-compact solutions for variants of the KdVB equation were constructed by Wazwaz\cite{wazwaz1}. Kudryashov showed that many of the solutions of the KdVB equation listed above can be converted to each other by using \cite{kudryashov1}. He also showed that the traveling wave solutions of the Fisher equation and the KdVB are in the same form and derived an exponential-type solution with the assistance of Weierstrass function for the KdVB equation.

\noindent
Besides proposed various exact solutions in the literature summarized above, some numerical techniques were also developed to solve some problems constructed with KdVB equation. One of the frontier numerical studies on KdVB equation is based on Bubnov-Galerkin's finite element method using cubic B-splines\cite{zaki1}. The temporal evaluation of a Maxwellian and the time evaluation of the solutions were studied in details in that study. A collocation method based on quintic B-spline functions were also derived for the numerical solutions of the KdVB equation\cite{zaki2}. A linear Galerkin-Fourier spectral technique was implemented for the numerical solutions of the KdVB eqution with periodic initial condition\cite{Lu1}.  A spectral collocation method based on differentiated Chebyshev polynomials combined with the fourth order Runge-Kutta method was proposed for the numerical solution of the initial boundary value problem modeling $\tanh$-type solitary wave solutions\cite{khaser1}. Unconditionally stable septic B-spline collocation method was proposed to solve the KdVB equation numerically by El-Danaf\cite{eldanaf1}.

\noindent
In this study, we consider some initial boundary value problems defined as
\begin{equation}
u_{t}+\varepsilon uu_{x}-\vartheta u_{xx}+\mu u_{xxx}=0 \ a< x< b
\label{kdvbivp}
\end{equation}
subject to the initial condition
\begin{equation*}
u(x,0)=f(x),\text{ }\  \  \ a\leq x\leq b
\end{equation*}
in the finite problem interval $[a,b]$. We choose the boundary conditions from the set
\begin{equation*}
\begin{aligned}
u(a,t)&=g_1(t), &u(b,t)&=g_2(t), \\ 
u_{x}(a,t)&=g_3(t), &u_{x}(b,t)&=g_4(t), \\ 
u_{xx}(a,t)&=g_5(t), &u_{xx}(b,t)&=g_6(t), \\ 
u_{xxx}(a,t)&=g_7(t), &u_{xxx}(b,t)&=g_8(t)
\end{aligned}
\end{equation*}
where $g_i(t),i=1,2,...,8$ denote $x$ independent functions.
\section{Method of Solution}
\subsection{Extended Cubic B-splines}
\begin{figure}[htbp]
	\centering
		\includegraphics[scale=0.5]{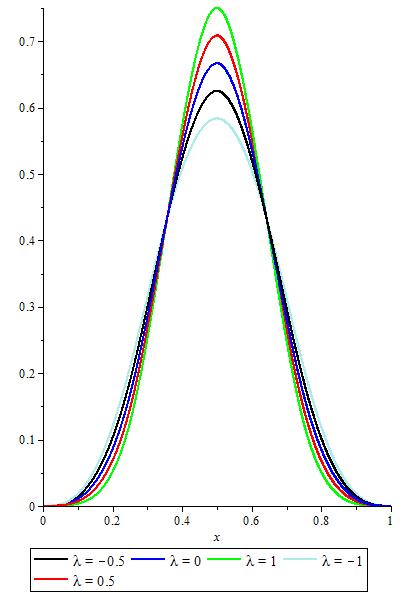}
	\caption{Extended B-splines for various values of the free parameter $\lambda$}
	\label{fig:Fig1}
\end{figure}

Let $\pi$ be a uniform grid distribution of the finite interval $[a,b]$ such as $\pi :a=x_{0}<x_{1}<\ldots <x_{N}=b$ with equal  sub interval length $h=(b-a)/N$. The extended form of a cubic B-spline $E_{i}$ is defined as \cite{prenter,dursunext}
\begin{equation}
E_{i}(x)=\frac{1}{24h^{4}}\left \{ 
\begin{array}{ll}
4h(1-\lambda )(x-x_{i-2})^{3}+3\lambda (x-x_{i-2})^{4}, & \left[
x_{i-2},x_{i-1}\right] , \\ 
\begin{array}{l}
(4-\lambda )h^{4}+12h^{3}(x-x_{i-1})+6h^{2}(2+\lambda )(x-x_{i-1})^{2} \\ 
-12h(x-x_{i-1})^{3}-3\lambda (x-x_{i-1})^{4}%
\end{array}
& \left[ x_{i-1},x_{i}\right] , \\ 
\begin{array}{l}
(4-\lambda )h^{4}-12h^{3}(x-x_{i+1})+6h^{2}(2+\lambda )(x-x_{i+1})^{2} \\ 
+12h(x-x_{i+1})^{3}-3\lambda (x-x_{i+1})^{4}%
\end{array}
& \left[ x_{i},x_{i+1}\right] , \\ 
4h(\lambda -1)(x-x_{i+2})^{3}+3\lambda (x-x_{i+2})^{4}, & \left[
x_{i+1},x_{i+2}\right] , \\ 
0 & \text{otherwise.}%
\end{array}%
\right.  \label{e1}
\end{equation}
where the real $\lambda$ is the free parameter. In fact, classical cubic polynomial B-splines are the particular form of extended cubic B-splines when $\lambda =0$. The set of the extended cubic B-splines $\{E_i(x) \}_{i=-1}^{N+1}$ is a basis for the functions defined in the interval $[x_0,x_N]$\cite{prenter,dursunext}. When the free parameter is chosen different from zero, the shape of the function changes slightly, Fig \ref{fig:Fig1}. The nonzero functional and derivative values of each extended cubic B-spline $E_{i}(x)$ at the grids of $\pi$in the $[a,b]$ can be summarized as in Table \ref{table1}.  

\begin{table}[htbp]
	\centering
	\caption{$E_i(x)$ and its derivatives at the grid points}
		\begin{tabular}{llllcl}
			\hline
$x$ & $x_{i-2}$ & $x_{i-1}$ & $x_{i}$ & $x_{i+1}$ & \multicolumn{1}{c}{$%
x_{i+2}$} \\ \hline
$24E_{i}(x)$ & $0$ & $4-\lambda $ & $16+2\lambda $ & $4-\lambda $ & 
\multicolumn{1}{c}{$0$} \\ 
$2hE_{i}^{^{\prime }}(x)$ & $0$ & $-1$ & $0$ & $1$ & \multicolumn{1}{c}{$0$} \\ 
$2h^{2}E_{i}^{^{\prime \prime }}(x)$ & $0$ & $2+\lambda $ & $-4-2\lambda $ & $%
2+\lambda $ & \multicolumn{1}{c}{$0$} \\ \hline
		\end{tabular}
	
	\label{table1}
\end{table}

\noindent
Having only first two derivatives of each extended B-spline $E_i(x)$ forces us to reduce the derivative order of the KdVB equation (\ref{kdvb}). Thus, defining a new variable $v(x,t)=u_x(x,t)$ reduces the KdVB equation (\ref{kdvb}) to a coupled system of equations
\begin{equation}
\begin{aligned}
u_{t}+\varepsilon uv-\vartheta v_{x}+\mu v_{xx}=0 \\ 
u_{x}-v=0 \label{sys}
\end{aligned}
\end{equation}
\subsection{Collocation Method}
\noindent
Let $U(x,t)$ and $V(x,t)$ be the approximate solutions to $u(x,t)$ and $v(x,t)$, respectively that
\begin{equation}
\begin{aligned}
U(x,t)&=\sum\limits_{i=-1}^{N+1}\delta_{i}E_{i}(x),\\
V(x,t)&=\sum _{i=-1}^{N+1}\phi_{i}E_{i}(x) 
\end{aligned} \label{e2}
\end{equation}%
where $\delta_{i}$ and $\phi_{i}$ are time dependent parameters that will be calculated via the extended cubic B-splines and the complementary conditions. The approximate solutions given in (\ref{e2}) and their first and second derivatives at a grid $x_{i}$  can be determined by using the Table (\ref{table1}). Thus, the approximate solutions and their derivatives take the form
\begin{equation}
\begin{aligned}
U_{i}&=U(x_{i},t)=\dfrac{4-\lambda }{24}\delta _{i-1}+\dfrac{8+\lambda }{12}\delta _{i}+\dfrac{4-\lambda }{24}\delta _{i+1}, \\ 
U_{i}^{\prime }&=U^{\prime }(x_{i},t)=\dfrac{-1}{2h}\left( \delta
_{i-1}-\delta _{i+1}\right) \\ 
U_{i}^{\prime \prime }&=U^{\prime \prime }(x_{i},t)=\dfrac{2+\lambda }{2h^{2}}\left( \delta _{i-1}-2\delta _{i}+\delta _{i+1}\right) \\
V_{i}&=V(x_{i},t)=\dfrac{4-\lambda }{24}\phi _{i-1}+\dfrac{8+\lambda }{12}\phi _{i}+\dfrac{4-\lambda }{24}\phi _{i+1}\\ 
V_{i}^{\prime }&=V^{\prime }(x_{i},t)=\dfrac{-1}{2h}\left( \phi _{i-1}-\phi_{i+1}\right)\\ 
V_{i}^{\prime \prime }&=V^{\prime \prime }(x_{i},t)=\dfrac{2+\lambda }{2h^{2}}\left( \phi _{i-1}-2\phi _{i}+\phi _{i+1}\right)
\label{e3}
\end{aligned}%
\end{equation}
To integrate the system (\ref{sys}) in time variable we first use forward finite difference and the Crank-Nicolson methods to give
\begin{equation}
\begin{array}{r}
\dfrac{U^{n+1}-U^{n}}{\Delta t}+\varepsilon \dfrac{(UV)^{n+1}+(UV)^{n}}{2}%
-\vartheta \dfrac{(V_{x})^{n+1}+(V_{x})^{n}}{2}+\mu \dfrac{%
(V_{xx})^{n+1}+(V_{xx})^{n}}{2}=0 \\ 
\dfrac{U_{x}^{n+1}+U_{x}^{n}}{2}-\dfrac{V^{n+1}+V^{n}}{2}=0%
\end{array}
\label{e5}
\end{equation}%
where $U^{n+1}=U(x,(n+1)\Delta t)$ represent the solution at the $(n+1)$.th
time level. Here $t^{n+1}$ $=t^{n}+t$, and $\Delta t$ is the time step,
superscripts denote $n$.th time level , $t^{n}=n\Delta t$. Linearizing the terms $(UV)^{n+1}$and $(UV)^{n}$ in (\ref{e5}) as
\begin{equation}
\begin{array}{l}
(UV)^{n+1}=U^{n+1}V^{n}+U^{n}V^{n+1}-U^{n}V^{n} \\ 
(UV)^{n}=U^{n}V^{n}%
\end{array}
\label{e6}
\end{equation}%
gives
\begin{equation}
\begin{array}{r}
\dfrac{U^{n+1}-U^{n}}{\Delta t}+\varepsilon (\dfrac{U^{n+1}V^{n}+U^{n}V^{n+1}%
}{2})-\vartheta (\dfrac{V_{x}^{n+1}+V_{x}^{n}}{2})+\mu (\dfrac{%
V_{xx}^{n+1}+V_{xx}^{n}}{2})=0 \\ 
\dfrac{U_{x}^{n+1}+U_{x}^{n}}{2}-\dfrac{V^{n+1}+V^{n}}{2}=0%
\end{array}
\label{e7}
\end{equation}%
Substitution of the approximate solutions (\ref{e2}) and their derivatives (\ref{e3}) into (\ref{e7}) yields
\begin{eqnarray}
&&\nu _{m1}\delta _{m-1}^{n+1}+\nu _{m2}\phi _{m-1}^{n+1}+\nu _{m3}\delta
_{m}^{n+1}+\nu _{m4}\phi _{m}^{n+1}+\nu _{m1}\delta _{m+1}^{n+1}+\nu
_{m5}\phi _{m+1}^{n+1}  \label{e9} \\
&=&\nu _{m6}\delta _{m-1}^{n}+\nu _{m7}\phi _{m-1}^{n}+\nu _{m8}\delta
_{m}^{n}+\nu _{m9}\phi _{m}^{n}+\nu _{m6}\delta _{m+1}^{n}+\nu _{m10}\phi
_{m+1}^{n}  \notag
\end{eqnarray}%
\begin{eqnarray}
&&\nu _{m11}\delta _{m-1}^{n+1}+\nu _{m12}\phi _{m-1}^{n+1}+0\delta
_{m}^{n+1}+\nu _{m13}\phi _{m}^{n+1}+\nu _{m11}\delta _{m+1}^{n+1}-\nu
_{m12}\phi _{m+1}^{n+1}  \label{e10} \\
&=&-\nu _{m11}\delta _{m-1}^{n}-\nu _{m12}\phi _{m-1}^{n}+0\delta
_{m}^{n}-\nu _{m13}\phi _{m}^{n}-\nu _{m11}\delta _{m+1}^{n}+\nu _{m12}\phi
_{m+1}^{n}  \notag
\end{eqnarray}
The coefficients in the system (\ref{e9}) and (\ref{e10}) are 
\begin{equation*}
\begin{array}{lll}
\nu _{m1}=\left( \dfrac{2}{\Delta t}+\varepsilon L\right) \alpha _{1} &  & 
\nu _{m8}=\dfrac{2}{\Delta t}\alpha _{2} \\ 
\nu _{m2}=\varepsilon K\alpha _{1}-\vartheta \beta _{1}+\mu \gamma _{1} &  & 
\nu _{m9}=-\mu \gamma _{2} \\ 
\nu _{m3}=\left( \dfrac{2}{\Delta t}+\varepsilon L\right) \alpha _{2} &  & 
\nu _{m10}=-\vartheta \beta _{1}-\mu \gamma _{1} \\ 
\nu _{m4}=\varepsilon K\alpha _{2}+\mu \gamma _{2} &  & \nu _{m11}=\beta _{1}
\\ 
\nu _{m5}=\varepsilon K\alpha _{1}+\vartheta \beta _{1}+\mu \gamma _{1} &  & 
\nu _{m12}=-\alpha _{1} \\ 
\nu _{m6}=\dfrac{2}{\Delta t}\alpha _{1} &  & \nu _{m13}=-\alpha _{2} \\ 
\nu _{m7}=\vartheta \beta _{1}-\mu \gamma _{1} &  & 
\end{array}%
\end{equation*}
where%
\begin{eqnarray*}
K &=&\alpha _{1}\delta _{i-1}+\alpha _{2}\delta _{i}+\alpha _{1}\delta _{i+1}
\\
L &=&\alpha _{1}\phi _{i-1}+\alpha _{2}\phi _{i}+\alpha _{1}\phi _{i+1}
\end{eqnarray*}
and%
\begin{eqnarray*}
\alpha _{1} &=&\dfrac{4-\lambda }{24},\text{ }\alpha _{2}=\dfrac{8+\lambda }{%
12} \\
\beta _{1} &=&-\dfrac{1}{2h},\text{ }\gamma _{1}=\dfrac{2+\lambda }{2h^{2}},%
\text{ }\gamma _{2}=-\dfrac{4+2\lambda }{2h^{2}}
\end{eqnarray*}%

The system (\ref{e9}-\ref{e10}) can be written in the following matrices system
\begin{equation}
\mathbf{Ax}^{n+1}=\mathbf{Bx}^{n}  \label{BST10}
\end{equation}%
where%
\begin{equation*}
\mathbf{A=}%
\begin{bmatrix}
\nu _{m1} & \nu _{m2} & \nu _{m3} & \nu _{m4} & \nu _{m1} & \nu _{m5} &  & 
&  &  \\ 
\nu _{m11} & \nu _{m12} & 0 & \nu _{m13} & \nu _{m11} & -\nu _{m12} &  &  & 
&  \\ 
&  & \nu _{m1} & \nu _{m2} & \nu _{m3} & \nu _{m4} & \nu _{m1} & \nu _{m5} & 
&  \\ 
&  & \nu _{m11} & \nu _{m12} & 0 & \nu _{m13} & \nu _{m11} & -\nu _{m12} & 
&  \\ 
&  &  & \ddots & \ddots & \ddots & \ddots & \ddots & \ddots &  \\ 
&  &  &  & \nu _{m1} & \nu _{m2} & \nu _{m3} & \nu _{m4} & \nu _{m1} & \nu
_{m5} \\ 
&  &  &  & \nu _{m11} & \nu _{m12} & 0 & \nu _{m13} & \nu _{m11} & -\nu
_{m12}%
\end{bmatrix}%
\end{equation*}%
and%
\begin{equation*}
\mathbf{B=}%
\begin{bmatrix}
\nu _{m6} & \nu _{m7} & \nu _{m8} & \nu _{m9} & \nu _{m6} & \nu _{m10} &  & 
&  &  \\ 
-\nu _{m11} & -\nu _{m12} & 0 & -\nu _{m13} & -\nu _{m11} & \nu _{m12} &  & 
&  &  \\ 
&  & \nu _{m6} & \nu _{m7} & \nu _{m8} & \nu _{m9} & \nu _{m6} & \nu _{m10}
&  &  \\ 
&  & -\nu _{m11} & -\nu _{m12} & 0 & -\nu _{m13} & -\nu _{m11} & \nu _{m12}
&  &  \\ 
&  &  & \ddots & \ddots & \ddots & \ddots & \ddots & \ddots &  \\ 
&  &  &  & \nu _{m6} & \nu _{m7} & \nu _{m8} & \nu _{m9} & \nu _{m6} & \nu
_{m10} \\ 
&  &  &  & -\nu _{m11} & -\nu _{m12} & 0 & -\nu _{m13} & -\nu _{m11} & \nu
_{m12}%
\end{bmatrix}%
\end{equation*}
There are $2N+2$ equations with $2N+6$ unknown parameters $$%
\mathbf{x}^{n+1}=(\delta _{-1}^{n+1},\phi _{-1}^{n+1},\delta _{0}^{n+1},\phi
_{0}^{n+1}\ldots ,\delta _{N+1}^{n+1},\phi _{N+1}^{n+1})$$ in this system. A unique solution of the system can be obtained by imposing the boundary conditions $U_{x}(a,t)=0,U_{x}(b,t)=0,V_{x}(a,t)=0,V_{x}(b,t)=0$ to have the following the equations:%
\begin{equation*}
\delta _{-1}=\delta _{1},\text{ }\phi _{-1}=\phi _{1},\text{ }\delta
_{N-1}=\delta _{N+1},\text{ }\phi _{N-1}=\phi _{N+1}
\end{equation*}
Elimination of the parameters $\delta _{-1},\phi _{-1},\delta _{N+1},\phi
_{N+1},$ using the equations (\ref{e9}) from the the system gives a solvable system of $2N+2$ linear equations including $2N+2$ unknown parameters. 

Since the right handside of this equation consists of only known values at the $n$.th time level, the solution of the system gives the solution at the $(n+1)$.th time level. In order to initialize this iteration system, we need the initial vector $\mathbf{x^0}$.The initial parameter vectors $d_{1}=(\delta
_{-1},\delta _{0},..\delta _{N},\delta _{N+1})$, $d_{2}=(\phi _{-1},\phi
_{0},..\phi _{N},\phi _{N+1})$ can be determined by using
$$
\begin{array}{l}
U_{xx}(a,0)=0=\gamma _{1}\delta _{-1}^{0}+\gamma _{2}\delta
_{0}^{0}+\gamma _{1}\delta _{1}^{0}, \\ 
U_{xx}(x_{i},0)=\gamma _{1}\delta _{i-1}^{0}+\gamma _{2}\delta
_{i}^{0}+\gamma _{1}\delta _{i+1}^{0}=U_{xx}(x_{i},0),i=1,...,N-1 \\ 
U_{xx}(b,0)=0=\gamma _{1}\delta _{N-1}^{0}+\gamma _{2}\delta
_{N}^{0}+\gamma _{1}\delta _{N+1}^{0}, \\ 
V_{x}(a,0)=0=\phi _{-1}^{0}-\phi _{1}^{0} \\ 
V_{x}(x_{i},0)=\phi _{i-1}^{0}-\phi
_{i+1}^{0}=V_{x}(x_{i},0),i=1,...,N-1 \\ 
V_{x}(b,0)=\phi _{N-1}^{0}-\phi _{N+1}^{0}%
\end{array}%
$$
\section{Numerical tests}
We chose some initial boundary value problems constructed on the the KdV-Burgers equation to check the accuracy and validity of the proposed method. The discrete maximum error norm
\begin{equation*}
L_{\infty }=\left \vert u-U\right \vert _{\infty }=\max \limits_{j}\left
\vert u_{j}-U_{j}^{n}\right \vert
\end{equation*}%
between the numerical and exact solution is computed for different values of the viscosity parameter $\vartheta $, different times with various time and space step sizes $\Delta t$, $h$. 
\subsection{Example 1: Initial boundary value problem with analytical solution}
Consider the initial boundary value problem constructed with the KdV-Burgers' equation (\ref{kdvb}) with the initial condition
\begin{equation}
u(x,0)=-\frac{6\vartheta ^{2}}{25\mu }\left[ 1+\tanh \left( \frac{\vartheta x}{10\mu } \right) +\frac{1}{2}\sech^{2}\left( \frac{\vartheta x}10\mu 
\right) \right]   \label{n1}
\end{equation}
and the boundary conditions
\begin{eqnarray*}
u(x,a) &=&-\frac{6\vartheta ^{2}}{25\mu }\left[ 1+\tanh \left( \frac{%
\vartheta }{10\mu }(a+\frac{6\vartheta ^{2}}{25\mu }t)\right) +\frac{1}{2}%
\sech^{2}\left( \frac{\vartheta }{10\mu }(a+\frac{6\vartheta ^{2}}{25\mu }%
t)\right) \right]  \\
u(x,b) &=&-\frac{6\vartheta ^{2}}{25\mu }\left[ 1+\tanh \left( \frac{%
\vartheta }{10\mu }(b+\frac{6\vartheta ^{2}}{25\mu }t)\right) +\frac{1}{2}%
\sech^{2}\left( \frac{\vartheta }{10\mu }(b+\frac{6\vartheta ^{2}}{25\mu }%
t)\right) \right] 
\end{eqnarray*}
These complementary conditions have been derived from the Fan and Zhang's \cite{fan1} analytical solution of the form 
\begin{equation}
u(x,t)=-\frac{6\vartheta ^{2}}{25\mu }\left[ 1+\tanh \left( \frac{\vartheta}{10\mu }\left( x+\frac{6\vartheta^2}{25 \mu} t \right) \right) +\frac{1}{2}\sech^{2}\left( \frac{\vartheta}{10\mu}\left( x+\frac{6\vartheta^2}{25 \mu} t \right)\right) \right]   \label{n1}
\end{equation}
for the particular choice of $\varepsilon =1$. This solution has three components, a constant, the $\tanh$ component and the $\sech$ component and it represents a traveling wave moving along the $x-$axis as time goes.

\begin{figure}[b]
    \subfigure[$\vartheta=0.001$]{
   \includegraphics[scale =0.5] {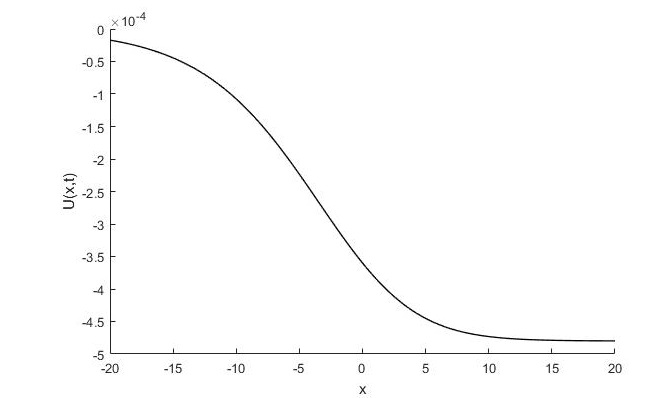}
   \label{fig:2a}
 }
 \subfigure[$\vartheta=0.005$]{
   \includegraphics[scale =0.5] {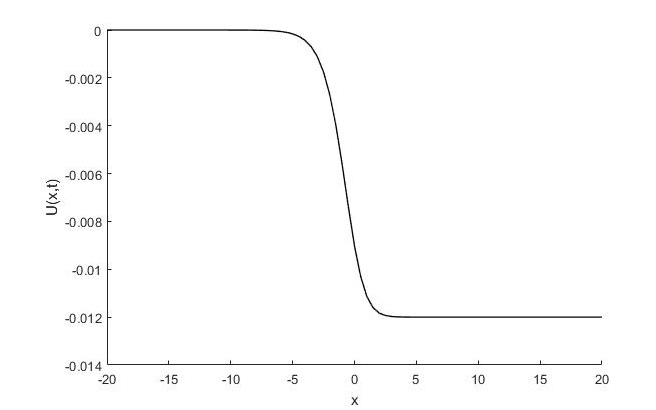}
   \label{fig:2b}
 }
  \subfigure[$\vartheta=0.01$]{
   \includegraphics[scale =0.5] {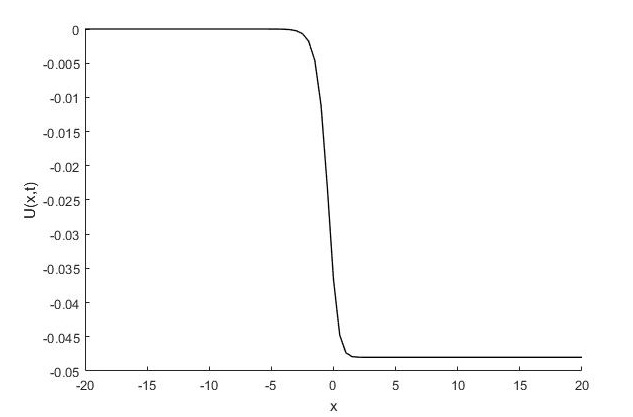}
   \label{fig:2c}
 }
 \caption{The profiles of travelling waves at $t=1$ for various values of $\vartheta$}
\end{figure}

\begin{figure}[b]
    \subfigure[$\vartheta=0.001$]{
   \includegraphics[scale =0.5] {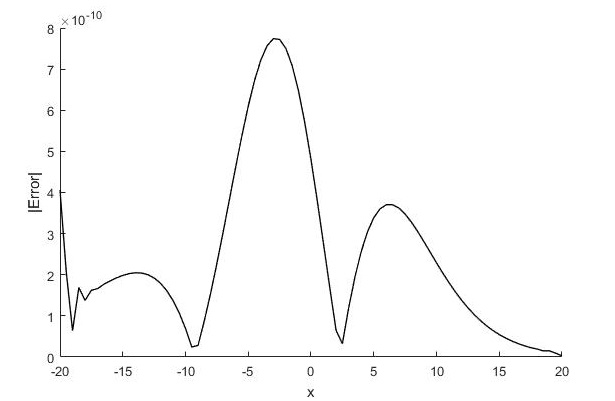}
   \label{fig:3a}
 }
 \subfigure[$\vartheta=0.005$]{
   \includegraphics[scale =0.5] {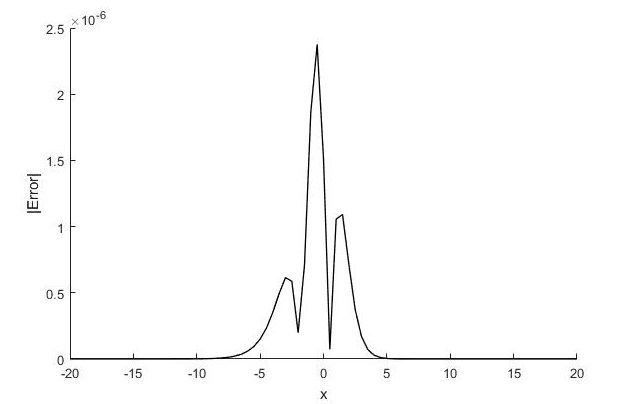}
   \label{fig:3b}
 }
  \subfigure[$\vartheta=0.01$]{
   \includegraphics[scale =0.5] {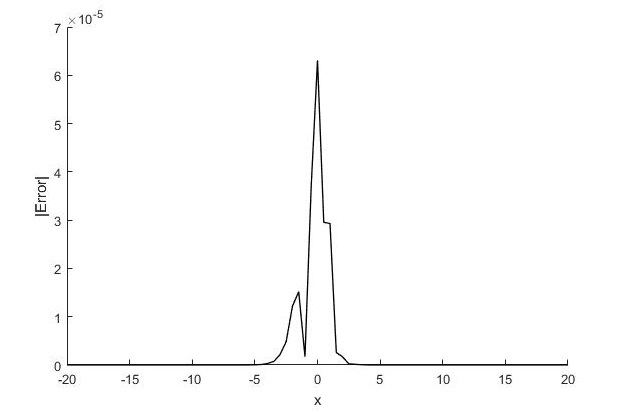}
   \label{fig:3c}
 }
 \caption{Error distributions at the time $t=1$ for various values of $\vartheta$}
\end{figure}

\noindent
For the sake of comparison with some earlier works, we fix the parameters $\mu =0.01$, $h=0.5$, $\Delta t=0.001$ and run the algorithm in the finite problem interval $[-20,20]$ up to the terminating time $t=1$. When the coefficient is chosen as $\vartheta =0.001, 0.005$ and $0.01$, the terminating time profiles of the waves are given in Fig \ref{fig:2a}-\ref{fig:2c}. The increase of the dispersion parameter value dominates the Burgers-type solution and the wave gets steeper. The maximum error distrubitions for each case are also plotted in Fig \ref{fig:3a}-\ref{fig:3c} for the same parameter values. The error distribution plots in all cases show that the error is greater near the descent part of the wave.

\noindent
A comparison of discrete maximum norms with the results obtained by the radial basis collocation methods with multiquadric (MQ), inverse quadric (IQ) and Gaussian (GA) forms in the study\cite{haq1} at the time $t=1$ and $t=10$ is tabulated in Table \ref{tab:test1}. For compatibility we chose the parameters as $\varepsilon=1$, $\vartheta =0.004$, $\Delta t=0.001$ and $h=0.5$. This comparison is an indicator of the efficiency of the proposed method based on the extended B-spline functions. Both the particular selections of the extention parameter of the extended B-splines as $\lambda =0$ and $\lambda =-1.969$ gives at least two decimal digits better results than the results of radial basis collocation methods. Even the proposed method for both values of extension parameter generate the results in the same decimal digit accuracy at the time $t=1$, the nonzero value of the extension parameter keeps the accuracy digits at $t=10$ but the classical cubic B-spline based method lose one decimal digit in accuracy. The method derived from radial basis functions can not catch the same accuracy even though the same parameters are used during the implementation of the proposed.

\begin{table}[b]
  \begin{center}
    \caption{Comparison of the maximum error with radial basis collocation methods}
    \label{tab:test1}
    \begin{tabular}{l|l|l|l|l|l}
    		\hline\hline
      	$t$ & Present ($\lambda =0$)& Present ($\lambda =-1.969$)& MQ \cite{haq1} & GA\cite{haq1} &IQ\cite{haq1} \\ \hline \hline
      	$1$ & $9.441 \times 10^{-11}$& $7.271 \times 10^{-11}$& $6.822 \times 10^{-9}$& $7.913 \times 10^{-9}$&  $4.077 \times 10^{-7}$ \\
      	$10$ & $1.269 \times 10^{-10}$& $9.509 \times 10^{-11}$& $2.479 \times 10^{-8}$& $3.294 \times 10^{-6}$&  $1.270 \times 10^{-6}$ \\
      \hline \hline
    \end{tabular}
  \end{center}
\end{table}

\noindent
\subsection{Example 2: Non analytical initial data}
In the second problem, we consider the particular form of the KdVB when $\vartheta =0$. Eliminating the dispersion term reduces the KdVB to the Korteweg-de Vries equation (KdVE). In order to check the validity of the method, we chose a non analytical solution describing the split of the initial pulse into a train of waves. This initial pulse is described by the data 
\begin{equation}
u(x,0)=\frac{1}{2}\left[ 1-\tanh{\frac{|x|-25}{5}}\right]
\end{equation}
and the routine is run up to the terminating time $t=800$ in the finite interval $[-50,150]$ for the compatibility to the results obtained in \cite{gardner1,korkmaz1,ozlem3}. During this long process time we will observe the lowest four conserved quantities \cite{miura1,korkmaz1}
\begin{equation*}
\begin{aligned}
C_1&=\int \limits_{-\infty}^{\infty}u dx\\ 
C_2&=\int \limits_{-\infty}^{\infty}u^{2}dx \\ 
C_3&=\int \limits_{-\infty}^{\infty}\left( u^{3}-\dfrac{3\mu }{\varepsilon }%
(u^{\prime })^{2}\right) dx \\
C_4&=\int \limits_{-\infty}^{\infty} \left( u^{4}-12\frac{\mu}{\varepsilon} u u_x^2+7.2\frac{\mu^2}{\varepsilon^2} u_{xx}^2\right)dx
\end{aligned}%
\end{equation*}
to be calculated for this initial boundary value problem. The approximate values of the conserved quantities are computed by modified trapezoid rule that the mean of the functional values at the consecutive points in the problem interval instead of the functional values at the grid points. The parameters are chosen as $\varepsilon =0.2$, $\mu =0.1$, $\Delta t=0.05$, $h=0.4$ in accordance with the earlier study of Korkmaz\cite{korkmaz1}. The conserved quantities for various choices of the extension parameter $\lambda$ are compared with the ones obtained by the differential quadrature methods based on cosine expansion (CDQ) and the Lagrange polynomials (LPDQ) combined with fourth order Runge-Kutta technique at various discrete times during the simulation, Table \ref{tab:con1}.
 
\noindent
It can be concluded that the lowest two quantities $C_1$ and $C_2$ that do not contain any derivatives of the solution are conserved successfully for all values of the extension parameter $\lambda$ even its zero. They both are in good agreement with the ones reported in Korkmaz' study \cite{korkmaz1}. Although the third lowest quantity $C_3$ is conserved during the simulation time for all choices of the extension parameter, its computed value changes dependently on the extension parameter $\lambda$. When the extension parameter approaches to $0$ in absolute value, its computed values also approache to its initial value. A similar situation is also observed during the calculation of $C_4$. Except the choice $\lambda =-1$, the other values of the extension parameter conserve the quantity acceptably even its numerical value deviates from its initial value. The deviations for all cases are directly proportional to the absolute value of extension parameter that the increase in the absolute value of the extension parameter enlarges the deviation. 
\begin{table}[b]
  \begin{center}
    \caption{Comparison of the lowest four conserved quantities}
    \label{tab:con1}{\scriptsize
    \begin{tabular}{l|l|l|l|l|l|r}
    		\hline\hline
    		Method&$\lambda$&$t$ & $C_1$ & $C_2$ & $C_3$ & $C_4$ \\
    		\hline \hline
    		 &&$0$&$50.000$&$45.000$&$42.301$&$40.442$\\
      \hline \hline
      Present&$0$&$200$&$50.001$&$45.002$&$42.404$&$40.983$\\
                &&$400$&$49.999$&$45.003$&$42.451$&$41.264$\\
                &&$600$&$49.999$&$45.003$&$42.453$&$41.290$\\
                &&$800$&$50.001$&$45.003$&$42.454$&$41.297$\\
      \hline
      Present&$-1$&$200$&$50.000$&$44.870$&$34.523$&$9.396$\\
                 &&$400$&$50.000$&$44.785$&$29.485$&$-7.116$\\
                 &&$600$&$50.000$&$44.773$&$28.803$&$-8.721$\\
                 &&$800$&$50.000$&$44.771$&$28.705$&$-8.904$\\
      \hline
      Present&$-0.5$&$200$&$50.000$&$44.957$&$39.748$&$30.237$\\
                    &&$400$&$50.000$&$44.928$&$38.050$&$24.631$\\
                    &&$600$&$50.000$&$44.924$&$37.811$&$24.061$\\
                    &&$800$&$50.000$&$44.923$&$37.776$&$23.991$\\
       \hline
      Present&$-0.25$&$200$&$50.000$&$44.983$&$41.263$&$36.353$\\
                    &&$400$&$50.000$&$44.971$&$40.557$&$34.078$\\
                    &&$600$&$50.000$&$44.969$&$40.455$&$33.846$\\
                    &&$800$&$50.000$&$44.969$&$40.439$&$33.814$\\
                           \hline
      Present&$-0.125$&$200$&$50.000$&$44.969$&$40.439$&$33.814$\\
                     &&$400$&$50.000$&$44.988$&$41.567$&$37.900$\\
                     &&$600$&$50.000$&$44.987$&$41.519$&$37.804$\\
                     &&$800$&$50.000$&$44.987$&$41.511$&$37.791$\\
                            \hline
      Present&$1$     &$200$&$50.000$&$45.047$&$45.072$&$51.883$\\
                     &&$400$&$49.997$&$45.078$&$46.889$&$58.291$\\
                     &&$600$&$50.001$&$45.082$&$47.140$&$58.933$\\
                     &&$800$&$50.004$&$45.082$&$47.174$&$58.994$\\    
                            \hline     
      Present&$0.5$   &$200$&$50.001$&$45.029$&$44.006$&$47.514$\\
                     &&$400$&$49.997$&$45.048$&$45.113$&$51.442$\\
                     &&$600$&$50.001$&$45.051$&$45.273$&$51.885$\\
                     &&$800$&$50.002$&$45.050$&$45.286$&$51.881$\\ 			
		                            \hline								
      Present&$0.25$  &$200$&$50.001$&$45.017$&$43.294$&$44.606$\\
                     &&$400$&$49.999$&$45.028$&$43.930$&$46.907$\\
                     &&$600$&$49.999$&$45.029$&$44.016$&$47.144$\\
                     &&$800$&$50.000$&$45.029$&$44.031$&$47.186$\\ 
		                            \hline								
      Present&$0.125$ &$200$&$50.001$&$45.001$&$42.786$&$42.899$\\
                     &&$400$&$49.998$&$45.016$&$43.233$&$44.248$\\
                     &&$600$&$50.001$&$45.017$&$43.280$&$44.387$\\
                     &&$800$&$50.000$&$45.017$&$43.286$&$44.396$\\ 
        \hline
     CDQ\cite{korkmaz1}&&$200$&$49.997$&$45.001$&$42.301$&$43.835$\\
                     &&$400$&$50.017$&$45.005$&$42.304$&$68.403$\\
                     &&$600$&$50.006$&$45.003$&$42.303$&$59.367$\\
                     &&$800$&$49.944$&$45.019$&$42.314$&$166.836$\\
                     \hline
      PDQ\cite{korkmaz1}&&$200$&$49.984$&$45.001$&$42.301$&$40.442$\\
                     &&$400$&$49.985$&$45.001$&$42.301$&$40.442$\\
                     &&$600$&$49.977$&$45.001$&$42.301$&$40.442$\\
                     &&$800$&$49.965$&$45.000$&$42.301$&$40.442$\\       \hline\hline
    \end{tabular}}
  \end{center}
\end{table}

\begin{figure}[htp]
    \subfigure[$t=0$]{
   \includegraphics[scale =0.5] {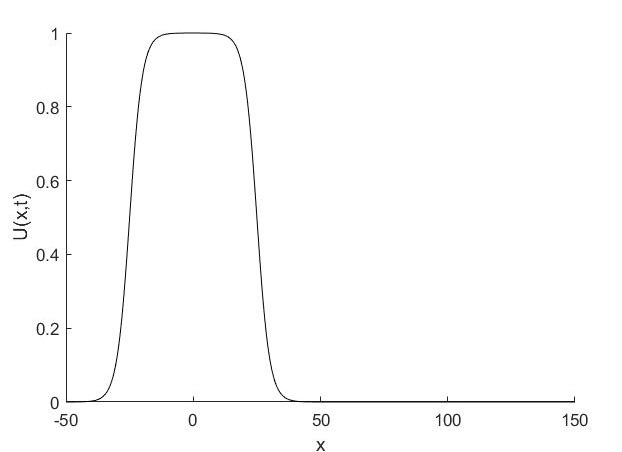}
   \label{fig:4}
 }
 \subfigure[$t=100$]{
   \includegraphics[scale =0.5] {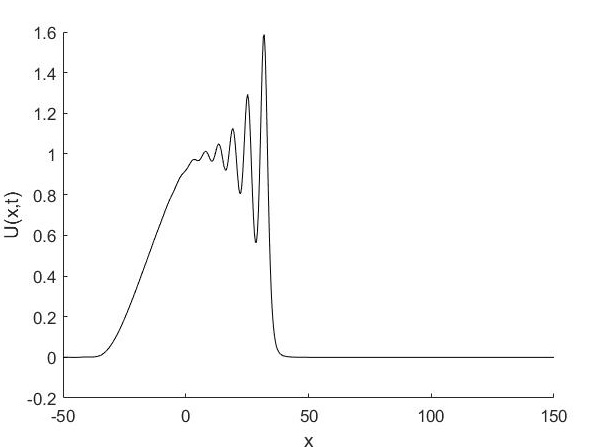}
   \label{fig:5}
 }
  \subfigure[$t=200$]{
   \includegraphics[scale =0.5] {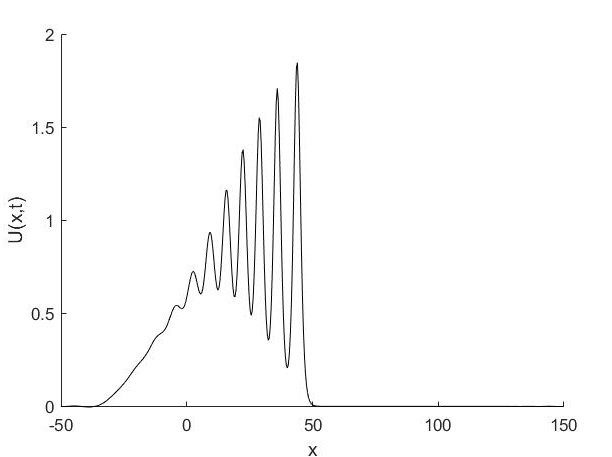}
   \label{fig:6}
 }
   \subfigure[$t=400$]{
   \includegraphics[scale =0.5] {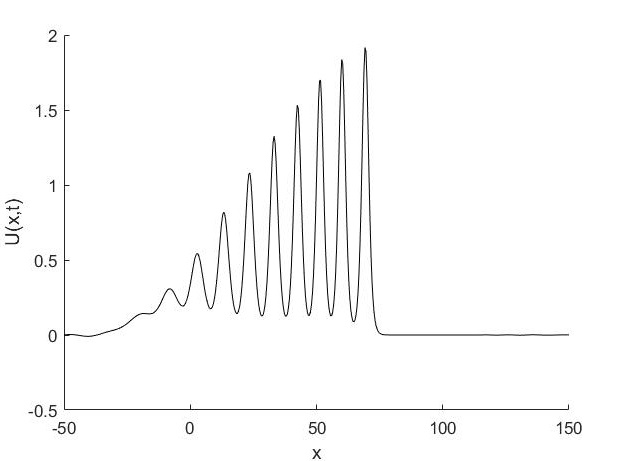}
   \label{fig:7}
 }
   \subfigure[$t=600$]{
   \includegraphics[scale =0.5] {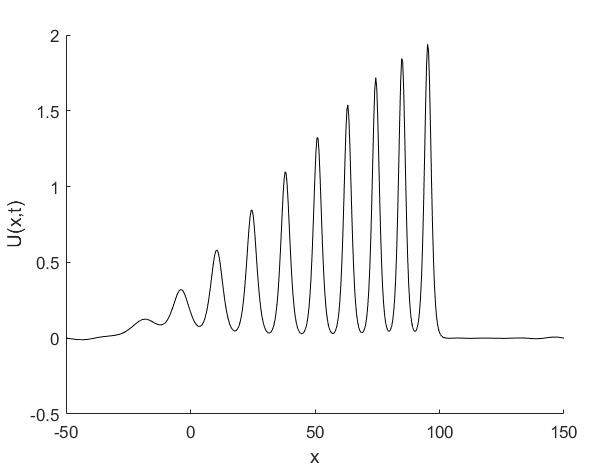}
   \label{fig:8}
 }
   \subfigure[$t=800$]{
   \includegraphics[scale =0.5] {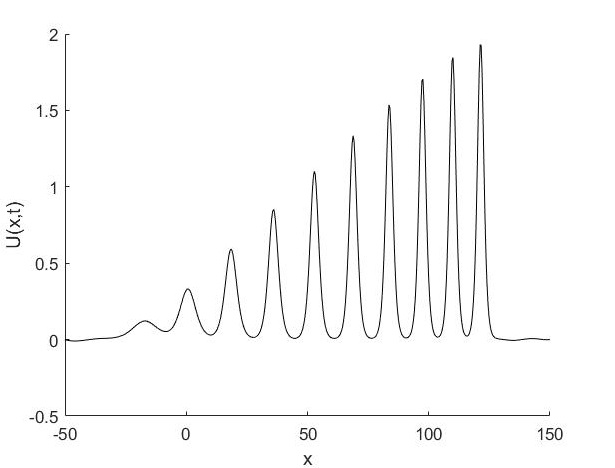}
   \label{fig:9}
 }
 \caption{The split of an initial pulse into waves}
\end{figure}

\noindent
The initial pulse initially positioned at $x=0$ begins its motion to the right as it splits into new waves. Three new waves and some oscillations following them can be observed at $t=100$. The leading wave's height is $1.587$ units and its peak is positioned at $x=32.0$. The peaks of the following two waves are positioned at $x=25.2$ and $19.2$ with heights $1.294$ and $1.126$, Fig \ref{fig:5}. The number of new apparent waves reaches $6$ at $t=200$. The peaks of the leading and the two followers of it reach $x=44.0$, $x=36.0$ and $x=28.8$, respectively. Their heights are measured as $1.846$, $1.701$ and $1.550$ at this distinct time, Fig \ref{fig:6}. At $t=400$, we see two more new waves are formed following the first $6$ waves, Fig \ref{fig:7}. When the time reaches $t=600$, the number of significant waves is $9$, Fig \ref{fig:8}. At the terminating time of the simulation, even  though we do not observe a significant change in the number of the formed waves, Fig \ref{fig:9}, an important increase in the velocity and heights of the leading waves is clearly seen. The first three leading waves reach $x=121.2$, $x=110.0$ and $x=97.60$ with heights $1.930$, $1.846$ and $1.703$, respectively.

\noindent

\section{Conclusion}

\noindent
The extended form of the classical polynomial cubic B-spline is adapted for the collocation method and implemented for the numerical solutions of some initial boundary value problems constructed with the KdVB equation even though the third order derivative of the extended cubic B-splines. The term $u_{xxx}$ is eliminated from the KdVB equation by a simple transformation. The coupled system derived from the elimination of $u_{xxx}$ is discretized in time by the Crank-Nicolson method. Following the linearization of the nonlinear term, the approximate solutions are substituted into the system instead of the exact solutions. After adapting the boundary conditions, the time iteration algorithm becomes ready to integrate the system in time variable.

\noindent 
The accuracy and efficiency of the proposed algorithm are monitored by solving two initial boundary value problems. The motion of traveling wave is considered in the first case. The solution simulations are accomplished successfully by the proposed method for various extension parameters $\lambda$. The accuracy of the results is measured by the discrete maximum norm that is showing the maximum distance between the analytical and numerical solutions. A comparison with some earlier studies show that the accuracy of the results obtained by the proposed method is at least two decimal digits better.

\noindent
In the second problem, we study wave formation from a non analytical initial data for the KdV equation obtained by eliminating the dispersion term. The graphical representations of the results at some distinct times are in a good agreement with some previous results. The lowest four conserved quantities are also computed at some specific times. The extension parameter has a critical role to increase the accuracy of the extended cubic B-spline based methods. 

\noindent
The optimum choice of the extension parameter in the cubic B-spline functions can be an interesting study for the future works. The required conditions to generate better results are the key for the improved numerical solutions.

\end{document}